\documentclass[10pt]{amsart}
\DeclareMathAlphabet{\mathscru}{OT1}{pzc}%
                                            {m}{it}\input mathrsfs.sty
\usepackage{verbatim,amsmath,latexsym,amssymb}

\theoremstyle{plain}
\newtheorem{thm}{Theorem}[section]

\newtheorem{lem}[thm]{Lemma}
\newtheorem{prop}[thm]{Proposition}
\newtheorem{cor}[thm]{Corollary}

\theoremstyle{definition}

\newtheorem{df}[thm]{Definition}

\newtheorem{ex}[thm]{Example}
\newtheorem{ex-notn}[thm]{Example/Notation}

\newtheorem{conj}[thm]{Conjecture}
\newtheorem{rem}[thm]{Remark}

\def\ann{\operatorname{ann}}

\def\conv{\operatorname{conv}}

\def\ext{\operatorname{Ext}}
\def\gr{\operatorname{gr}}
\def\hom{\operatorname{Hom}}

\def\id{\operatorname{id}}
\def\im{\operatorname{im}}

\def\mods{\operatorname{-mods}}
\def\rk{\operatorname{rk}}
\def\ord{\operatorname{ord}}

\def\qdeg{\operatorname{qdeg}}

\def\shom{\operatorname{{\mathscr H\hskip-0.6ex\mathit{om}}}}

\def\tdeg{\operatorname{tdeg}}

\def\var{\operatorname{Var}}
\def\vol{\operatorname{vol}}
 
\def\eps{\varepsilon}
\def\del{\partial}

\def\bolda{{\boldsymbol a}}
\def\bolde{{\boldsymbol e}}

\def\boldu{{\boldsymbol u}}
\def\boldv{{\boldsymbol v}}
\def\boldw{{\boldsymbol w}}

\def\calD{{\mathcal D}}
\def\calH{{\mathcal H}}

\def\calK{{\mathcal K}}

\def\calO{{\mathcal O}}

\def\fraka{{\mathfrak a}}
\def\frakm{{\mathfrak m}}

\def\CC{{\mathbb C}}    
\def\DD{{\mathbb D}}

\def\NN{{\mathbb N}}

\def\RR{{\mathbb R}}
\def\ZZ{{\mathbb Z}}

\def\bar#1{\overline{#1}}

\def\into{\hookrightarrow}
\def\Mtwo{{\em Macaulay} 2\expandafter}
\def\onto{\to\hskip-1.7ex\to}
\def\st{\operatorname{|}}
    
\def\ideal#1{{\langle #1 \rangle}}

\def\mylabel#1{\label{#1}}
\def\ignore#1{}

\def\comment#1{}
\numberwithin{equation}{section}
\begin{document}
 
\title[Duality and reducibility]
       {Duality and monodromy reducibility of $A$-hypergeometric systems}
\author[Uli Walther]{Uli Walther}
\address{Purdue University and Universit\"at Leipzig}
\email{walther@math.purdue.edu}
\begin{abstract}
We study $A$-hypergeometric systems $H_A(\beta)$ in the sense of
Gelfand, Kapranov and Zelevinsky under two aspects: the structure of
their holonomically dual system, and reducibility of their rank
module. 

We prove first that rank-jumping parameters always correspond to
reducible systems, and we show that the property of being reducible is
``invariant modulo the lattice''. In the second part we study a conjecture
of Nobuki Takayama which states that the holonomic dual of $H_A(\beta)$ is
of the form $H_A(\beta')$ for suitable $\beta'$.  We prove the
conjecture for all matrices $A$ and generic parameter $\beta$, exhibit
an example that shows that in general the conjecture cannot hold, and
present a refined version of the conjecture.

Questions on both duality and reducibility have been impossible to
answer with classical methods. This paper may be seen as an example of
the usefulness, and scope of applications, of the homological tools
for $A$-hypergeometric systems developed in \cite{MMW}.

\end{abstract}
\thanks{The author wishes to thank the NSF, the DfG and the Humboldt
foundation for their support.}
 
\maketitle

\section{Introduction}

Let $A$ be a $d\times n$-matrix with integer entries and distinct
columns $\{\bolda_j\}$; we abuse notation and denote by $A$ also the set of
its columns.  We assume throughout that $\ZZ A=\ZZ^d$, and that $A$ is
pointed: the semigroup $\NN A$ is to contain to units besides the
origin.  We denote by $x_A=\{x_j\st\bolda_j\in A\}$ a set of $n$
indeterminates, and by $\del_A=\{\del_j\st\bolda_j\in A\}$ the corresponding
partial differentiation operators.  Let $D_A$ be the Weyl algebra
$\CC\ideal{x_A,\del_A }$ and by $R_A$ the polynomial subring
$\CC[\del_A]$.

Gelfand, Kapranov and Zelevinsky discuss in \cite{GKZ-systems} the following
$D_A$-module $H_A(\beta)$. Let $I_A$ be the kernel of the
map $R_A\to \CC[t_1,\ldots,t_d]$ sending $\del_j$ to $t^{\bolda_j}$ where
here and henceforth we freely use multi-index notation. For $1\le i\le
d$ put $E_i=\sum_{j=1}^n a_{i,j}x_j\del_j$ and let
$\beta\in\CC^d$. Then by \cite{GKZ-systems} the {\em
$A$-hypergeometric system} (or GKZ-system) is
\[
H_A(\beta)=D_A/D_A\cdot(I_A,\{E_i-\beta_i\}_{i=1}^d).
\]
This module is holonomic, and it has regular singularities if 
$I_A$ defines a projective variety (i.e., if $I_A$ is homogeneous
in the usual sense). 

It is known by \cite{SST} that the holonomic rank of $H_A(\beta)$, the
$\CC(x_A)$-dimension of
$H_A(\beta)(x):=\CC(x_A)\otimes_{\CC[x_A]}H_A(\beta)$, is bounded from
below by the simplicial volume of $A$, the quotient of the Euclidean
volume of the convex hull of $A$ and the origin by the Euclidean
volume of the standard simplex in $\ZZ A=\ZZ^d$. This bound is sharp
for all parameters in a Zariski open set of $\CC^d$, and the structure
of the exceptional, {\em rank-jumping}, parameters are described in
\cite{MMW}.

In this paper we study $H_A(\beta)$ from two aspects of $D$-module
theory: duality and reducibility.

\subsection{Duality}
On $D_A$ there is a standard transposition which in multi-index
notation reads $\tau: x^\boldu\del^{\boldv}\mapsto
(-\del)^{\boldv}x^\boldu$; it induces an equivalence between
the categories of left and right $D_A$-modules.  The {\em holonomic
dual} of a holonomic left $D_A$-module $M$ is by definition
\[
\DD(M)=\tau(\ext^n_{D_A}(M,D_A)).
\]
Let $\calO^{an}_{\CC^n}$ be the analytic structure sheaf on $\CC^n$
and $\calD^{an}_{\CC^n}$ its sheaf of $\CC$-differential operators. On
the level of (constructible) solution sheaves
$\RR\shom_{\calD^{an}_{\CC^n}}(\calD^{an}_{\CC^n}\otimes_{D_A}
M,\calO^{an}_{\CC^n})$ the holonomic dual corresponds to Verdier
duality.

N.~Takayama conjectured that the holonomic dual of a GKZ-system
$H_A(\beta)$ is another GKZ-system. We discuss this conjecture and
prove that for any $A$ and for generic parameters Takayama's
conjecture holds. More precisely, we prove:
\begin{itemize}
\item If $\NN A$ is a normal semigroup then the set of parameters
  $\beta$ for which $\DD(H_A(\beta))$ fails to be a GKZ-system is
  contained in a finite subspace arrangement of codimension at least three
  with rational coefficients. 
\item For arbitrary $A$, the set of parameters $\beta$ for which
  $\DD(H_A(\beta))$ is not a GKZ-system is contained in a finite
  hyperplane arrangement with rational coefficients.
\item If $\dim(S_A)=2$ and $\beta$ is $A$-exceptional in the sense
  that $\rk(H_A(\beta))$ is larger than usual (see \cite{MMW} and
  below for details) then $\DD(H_A(\beta))$ is not a GKZ-system.
\end{itemize}

\subsection{Reducibility}
A $D_A$-module $M$ is said to have {\em irreducible monodromy
representation} if the {\em rank module}
$M(x):=\CC(x_A)\otimes_{\CC[x_A]}M$ is an irreducible module over
$D_A(x):=(\CC(x_A)\otimes_{\CC[ x_A]}D_A)$.  Gelfand, Kapranov and
Zelevinsky proved in \cite{GKZ-red} that an $A$-hyper\-geometric
system $H_A(\beta)$ has irreducible monodromy representation if the
parameter $\beta$ is generic.  In this note we prove that if $\beta$
has the property that $H_A(\beta)$ has irreducible monodromy
representation then all $\gamma\in \beta+\ZZ^d$ have the same
property; hence irreducibility descends to the quotient of $\CC^d$
modulo $\ZZ A$. In particular, all $\beta\in\ZZ^d$ give reducible
systems. The basic principle of proof for these facts are the {\em
contiguity operators} $\del_j:H_A(\beta)\to H_A(\beta+\bolda_j)$ used
in conjunction with the mechanisms of Euler--Koszul homology from
\cite{MMW}.

We further show with similar methods that
if $\beta$ is a rank-jumping parameter for $A$ then $H_A(\beta)$ must
have reducible monodromy.  We show this by producing a natural
morphism from $H_A(\beta)$ whose kernel exhibits a nontrivial submodule of
$\CC(x_A)\otimes_{\CC[x_A]} H_A(\beta)$.

\section{Toric modules and Euler--Koszul homology}

In this section we recall some notions and results from \cite{MMW}
regarding toric modules, Euler--Koszul homology, and generalized
$A$-hypergeometric systems. For notation, let us write $\calO_A$ for
$\CC[x_j\st\bolda_j\in A]$. 
\begin{df}
We $\ZZ A$-grade $R_A$ and $D_A$ according to $x_j\to \bolda_j$, $\del_j\to
-\bolda_j$.  The $i$-th {\em degree} function we denote by $\deg_i(-)$ so that
$\deg(-)=(\deg_1(-),\ldots,\deg_d(-))$ and the degree of $\del^\boldu\in
R_A$ is $A\cdot \boldu$. Note that
\[
I_A=\ideal{\del^\boldu-\del^\boldv\st \boldu,\boldv\in\NN^n, A\cdot
  \boldu=A\cdot \boldv}. 
\]
The {\em canonical degree} is $\eps_A=\sum
\bolda_j$, the degree of the generator of the canonical module of $R_A$. We
shall use ``graded'' and ``homogeneous'' to mean ``$\ZZ A$-graded''
and ``homogeneous with respect to the $\ZZ A$-grading''. We note that
we do not require $(1,1,\ldots,1)$ to be in the rowspan of $A$, as
$I_A$ may well not be projective. 

Let $S_A$ be the toric ring $R_A/I_A$; since $I_A$ is homogeneous, it
inherits the $\ZZ A$-grading from $R_A$.  A graded $R_A$-module $M$ is
called a {\em toric module} if it has a finite filtration by $A$-graded
$R_A$-modules such that each filtration quotient is a finitely
generated $S_A$-module. 
We say that $\beta\in\CC^d$ is a
{\em true degree} of the $A$-module $M$ if $\beta$ is contained in the
set of degrees of $M$:
\[
\left[ \beta\in\tdeg(M)\right] \Longleftrightarrow 
\left[\beta\in {\deg(m)\st 0\not =m\in
  M} \right].
\]
Moreover, we let the {\em quasi-degrees} of $M$ be the points in the
Zariski closure of $\tdeg(M)$.  For every toric module the
quasi-degrees form a finite subspace arrangement where each
participating subspace is a shift of a face of $\RR_{\geq 0}A$ by a
lattice element.

The toric
modules with $\ZZ A$-homogeneous maps of degree zero form the category
$A\mods$ which is closed under graded subquotients and extensions.
\end{df}

For any toric $R_A$-module $M$ one can define a collection of $d$
commuting endomorphisms denoted $E_i-\beta_i$, $1\le i\le d$, on
$D_A\otimes_{R_A}M$ which operate on an $A$-homogeneous element $m\in
D_A\otimes_{R_A} M$ by
$m\mapsto (E_i-\beta_i)\circ m$, where
\[
(E_i-\beta_i)\circ m=
(E_i-\beta_i-\deg_i(m))m.
\]
There is a functor $\calK_\bullet(\beta;-)=\calK_\bullet(E-\beta;-)$
from the category of $A$-graded $R_A$-modules to the category of
complexes of $\ZZ A$-graded $D_A$-modules; its output are complexes
concentrated between homological degrees $0$ and $d$.  It sends $M$ to
the Koszul complex defined by all morphisms $E_i-\beta_i$.  The
functor restricts to the category of toric modules on which it returns
complexes with holonomic homology.  By \cite{MMW}, if $-\beta$ is not
a quasi-degree of $M$ then $\calK_\bullet(\beta;M)$ is exact, while
if $\dim(M)=d$ and $M$ is Cohen--Macaulay then
$\calK_\bullet(\beta;M)$ is a resolution of its $0$-th homology
module.  A short exact sequence
\[
0\to M'\to M\to M''\to 0
\]
in $A\mods$  induces a short exact sequence of Koszul complexes
\[
0\to \calK_\bullet(\beta;M')\to \calK_\bullet(\beta;M)\to
\calK_\bullet(\beta;M'')\to 0
\]
which in turn produces a long exact sequence of {\em Euler--Koszul
homology}
\[
\cdots\to \calH_i(\beta;M'')\to \calH_{i-1}(\beta;M')\to
\calH_{i-1}(\beta;M)\to \calH_{i-1}(\beta;M'')\to\cdots
\]
where we have put $\calH_i(\beta;-)=H_i(\calK_\bullet(E-\beta;-))$. A
module of the form $\calH_0(\beta;M)$ is called a {\em generalized
GKZ-system}; if $M=S_A$ then $\calH_0(\beta;M)$ is a {\em proper} (or
{\em classical}\,) GKZ-system.

The grading interacts with Euler--Koszul homology as follows: for all
$\alpha\in\ZZ A$, 
\[
\calH_i(\beta;M(\alpha))=\calH_i(\beta-\alpha;M)(\alpha)
\]
where the left hand side is Euler--Koszul homology of the shifted
module $M$ while the right hand side is the shifted Euler--Koszul
homology of $M$ with shifted parameter. 

Finally, for all toric modules $M$ there is a {\em duality spectral
sequence}
\begin{eqnarray}
\mylabel{eqn-ss} \calH_i(-\beta-\eps_A;\ext^j_{R_A}(M,R_A))\Longrightarrow
\DD(\calH_{j-i}(\beta;M))^\vee,
\end{eqnarray}
where $M^\vee $ is the result of applying the coordinate change
$x\to -x$ to $M$.

\section{Reducibility of hypergeometric systems}

For any $\CC[x]$-module $H$ let $H(x)$ denote
$\CC(x)\otimes_{\CC[x]} H$ and recall that $(-)(x)$ is an exact
functor. We use similar notation for morphisms, if $\phi:M\to N$ then
$\phi(x):M(x)\to N(x)$.  In this section we investigate when the module
$H_A(\beta)(x)$ has nontrivial submodules; in the absence of such
submodules one says that $H_A(\beta)$ has {\em irreducible
monodromy}. If $H_A(\beta)(x)$ is irreducible and if
$\phi\in\hom_{D_A}(H_A(\beta),\CC\{x-c\})$ is any solution to
$H_A(\beta)$ at the point $c\in \CC^n$, then the $D_A$-annihilator of
$\phi(1)$ annihilates every single solution of $H_A(\beta)$. It is in
this sense that the monodromy is irreducible: analytic continuation of
any solution germ discovers all the others.

In our study of irreducibility we shall encounter the contiguity
operators $\del_k:H_A(\beta-\bolda_k)\to H_A(\beta)$.
 
Recall that there is a $(0,1)$-filtration on $D_A$ that leads to an
associated graded ring isomorphic to a ring of polynomials in $2n$
variables where the image of each $x_j$ has degree zero, and that of
each $\del_j$ has degree one. This filtration relates to holonomic
rank by the formula (see \cite{SST})
\[
\rk(H)=\dim_{\CC(x)}(\gr_{(0,1)}(H)(x))=\dim_{\CC(x)}(H(x)) 
\]
for any holonomic module $H$.

\subsection{Reducibility and rank jumps}

In this subsection we show that if a GKZ-system is rank-jumping then
it cannot have irreducible monodromy. We do this by showing that the
natural inclusion of $S_A$ into its normalization $\tilde S_A$ exhibits a
nontrivial submodule of $H_A(\beta)(x)$.

\begin{lem}
Consider the generalized hypergeometric system $\calH_0(\beta;\tilde
S_A)$. Its rank is equal to $\vol(A)$, independently of
$\beta$.
\end{lem}
\begin{proof}
We use the short exact sequence
\[
0\to S_A\to \tilde S_A\to C\to 0
\]
where $C$ is the $\ZZ A$-graded cokernel of the natural embedding of
$S_A$ into $\tilde S_A$. Normal seimgroup rings being Cohen--Macaulay, we
obtain a natural 4-term exact sequence of Euler--Koszul homology
\begin{eqnarray}
\label{eqn-4-term-C}
0\to \calH_1(\beta;C)\to \calH_0(\beta;S_A)\stackrel{\psi}{\longrightarrow}
 \calH_0(\beta;\tilde S_A)\to \calH_0(\beta;C)\to 0.
\end{eqnarray}
Since $S_A$ contains an ideal isomorphic to a graded shift of $\tilde
S_A$, $\dim(C)<\dim(S_A)$. In particular, for generic parameters
$\beta$ one has 
$-\beta\not \in \qdeg(C)$, the outer terms of the 4-term
sequence are zero and the map $\psi$ is an isomorphism. This means
that for generic parameters, $\rk(\calH_0(\beta;\tilde
S_A))=\rk(\calH_0(\beta;S_A))=\vol(A)$.

Now $\tilde S_A$ is a Cohen--Macaulay ring and hence a Cohen--Macaulay
$S_A$-module. By Corollary 9.2 in \cite{MMW}, $\tilde S_A$ has no
exceptional parameters. Hence all $\beta$ give the same
$\rk(H_A(\beta))$, and as generic ones give rank $\vol(A)$, all
parameters must.
\end{proof}

\begin{rem}
Pick a set of generators $\tilde A$ for the normalization of $\NN A$.
The volume of $A$ is usually not equal to the volume of $\tilde A$ --
even if $\tilde A$ is chosen minimally. So the expression
$\rk(\calH_0(\beta;\tilde S_A))$ depends strongly on the chosen base
ring $R$. Unless indicated otherwise, we consider $R_A$-modules.
\end{rem}

Now let $\beta$ be rank-jumping for $A$, $\rk(H_A(\beta))>\vol(A)$.
Since $\tilde S_A$ is Cohen--Macaulay, $\calH_1(\beta;\tilde
S_A)=0$. With $C=\tilde S_A/S_A$, there is an exact sequence
of type (\ref{eqn-4-term-C}).
By our lemma and the rank-jumping assumption,
$\ker(\phi)(x)\not =0$. We shall show that $\im(\phi)(x)\not=0$ and thus
prove
\begin{thm}
\mylabel{thm-jump-red}
If $\beta$ is rank-jumping for $A$ then $H_A(\beta)$ has reducible
monodromy; that is,
$H_A(\beta)(x)$ has a nontrivial composition chain.
\end{thm}
The proof of this theorem will occupy the remainder of this subsection.

\begin{proof}
Suppose the theorem is false and let $(A,\beta)$ be a counterexample
with minimal number of columns. Then $\rk(H_A(\beta))>\vol(A)$ and
$H_A(\beta)(x)$ is irreducible.  In the light of the fact that
$\ker(\phi)(x)\not =0$, that means that
$\calH_1(\beta;C)(x)=\calH_0(\beta;S_A)(x)$.

Since $\dim(C)<d$, $C$ has a composition chain whose factors are
isomorphic to toric rings associated to proper faces of $A$. In
particular, the $S_A$-module $C$ is annihilated by some monomial
$\del^\boldu\in R_A$. That implies that the submodule of
$\calH_1(\beta;C)=\calH_0(\beta;S_A)(x)$ generated by $\del^\boldu$ is
the zero module.

\begin{df}
For a fixed matrix $A$ and a parameter $\beta$ we set
$\fraka_{\beta}$ to be the subset of  $R_A$ given by
\begin{eqnarray*}
\left[ P\in \fraka_{\beta}\right]
&\Leftrightarrow &
\left[P\in R_A\st P=0 \text{ in } \calH_0(\beta;S_A)(x)\right]\\
&\Leftrightarrow &
[P\in R_A\st \exists f\in\calO_A, fP=0\in H_A(\beta)]\\
&\Leftrightarrow &
[P\in R_A\st \exists f\in\calO_A, fP\in D_A\cdot(I_A,\{E_i-\beta_i\}_{i=1}^d)]
\end{eqnarray*}
and note that $\fraka_\beta\supseteq I_A$ for all $A,\beta$. 
\end{df}
We now study general properties of $\fraka_{\beta}$. 
\begin{lem}
The set $\fraka_{\beta}$ is an ideal of $R_A$. 
\end{lem}
\begin{proof}
Obviously, $\fraka_\beta$ is a $\CC$-vector space. Let
$P\in\fraka_{\beta}$ with $fP\in D_A(I_A,\{E_i-\beta_i\}_{i=1}^n)$, 
and $Q\in D_A$. We shall prove that
$f^{\ord(Q)+1}QP=0$ in $\calH_0(\beta;S_A)$, which in particular implies
the lemma.

The claim is clearly true if the order of $Q$ is zero since then $f$
and $Q$ commute. Assume that $Q$ is order $k>0$ and the claim has been
proved for operators $Q$ up to order $k-1$. We have
\[
fQP=[f,Q]P+QfP=[f,Q]P
\]
in $\calH_0(\beta;S_A)$. Now the operator $[f,Q]$ has order smaller
than $k$. It follows then by induction on the order of $Q$ that
\[
0=f^{k-1+1}[f,Q]P=f^{k+1}QP
\]
in $\calH_0(\beta;S_A)$.  In particular, if $Q\in R_A$ then $QP\in
\fraka_\beta$.
\end{proof}

If $P\in\fraka_{\beta}$ then we can write it as a finite sum 
$P=\sum_{\alpha\in\ZZ A} P_\alpha$ where
$P_\alpha$ is homogeneous of degree $\alpha$. Then since
$P$ and $E_i-\beta_i$ are zero in $H_A(\beta)(x)$ we get that
\[
H_A(\beta)(x)\ni 0=(E_i-\beta_i) P-P(E_i-\beta_i)=[E_i,P]=\sum_{\alpha\in\ZZ A}
\alpha_i P_\alpha.
\]
In particular, $\sum_{\alpha\in\ZZ A} \alpha_i
P_\alpha\in\fraka_\beta$.  Repeated iteration implies that if
$P\in\fraka_{\beta}$, then $\sum_{\alpha\in\ZZ A} {\alpha_i}^k P_\alpha\in
\fraka_{\beta}$.  The nonvanishing of the appropriate Vandermonde
determinant assures that for all integers $a$ the subsum 
$\sum_{\alpha\in\ZZ A\st \alpha_i=a}
P_\alpha$ is in $\fraka_{\beta}$. Repeating this argument
for all $i$ it follows that $P_\alpha\in\fraka_\beta$. In particular,
$\fraka_{\beta}$ is a $\ZZ A$-graded ideal containing $I_A$.

Since in $S_A=R_A/I_A$ homogeneous elements of the same degree are constant
multiples of each other, $\fraka_{\beta}$ is an ideal of $R_A$
generated by $I_A$ and monomials.

\begin{prop}
\mylabel{prop-prime-mu}
Suppose $H_A(\beta)(x)$ is an irreducible module.  Then there exists a
monomial $\del^\boldu$ with $A\cdot \boldu=\alpha$ such that
$\fraka_{\beta-\alpha}$ is a prime ideal containing $\fraka_{\beta}$.
\end{prop}
\begin{proof}
Since $I_A$ is prime, we may assume that $\fraka_{\beta}$ strictly
contains $I_A$
and that there is a product
$\del_k\cdot\del^{\boldv-\bolde_k}=\del^\boldv\in \fraka_{\beta}$
such that neither $\del_k$ nor $\del^{\boldv-\bolde_k}$ is in
$\fraka_\beta$.  There is a sequence
\[
0\to S_A(\bolda_k)\to S_A\stackrel{\del_k}{\longrightarrow}
S_A/\ideal{\del_k}\to 0.
\]
We study the
associated long exact sequence of Euler--Koszul homology:
\[
\calH_1(\beta;S_A/\ideal{\del_k})\to H_A(\beta-\bolda_k)
\stackrel{\phi}{\longrightarrow} 
H_A(\beta)\to
\calH_0(\beta;S_A/\ideal{\del_k})\to 0.
\]
The image of $\phi(x)$ is generated by the coset of $\del_k$ in
$H_A(\beta)(x)$. Since $\del_k$ is not in $\fraka_\beta$, $\phi(x)$ is
not the zero map and since furthermore $H_A(\beta)(x)$ is irreducible,
$\phi(x)$ is onto.  This in turn implies that
$\calH_0(\beta;S_A/\ideal{\del_k})$ has zero rank. By Proposition 5.3
in \cite{MMW}, $\calH_i(\beta;S_A/\ideal{\del_k})=0$ for all
$i$. Hence $\phi$ and $\phi(x)$ are isomorphisms. Since
$\phi$ is right multiplication by $\del_k$, any monomial $\del^\boldw$
is in $\fraka_{\beta-\bolda_k}$ if and only if $\del^\boldw\cdot\del_k$ is
in $\fraka_{\beta}$. In other words,
\begin{eqnarray}\label{eqn-quotient}
\fraka_{\beta-\bolda_k}&=&
\fraka_{\beta}\st \ideal{\del_k}\supseteq \fraka_\beta
\end{eqnarray}
and by construction we have $\del^{\boldv-\bolde_k}\in
\fraka_{\beta-\bolda_k}\setminus \fraka_\beta$.

Iterating this procedure with $\beta-\bolda_k$ we get an increasing
sequence of ideals in $R_A$. Observing that $R_A$ is Noetherian we
arrive eventually at a parameter $\beta-\alpha\in\beta-\NN A$ such
that for every
$\del^{\boldv-\bolde_k}\cdot\del_k=\del^\boldv\in\fraka_{\beta-\alpha}$
either $\del_k\in \fraka_{\beta-\alpha}$ or
$\del^{\boldv-\bolde_k}\in\fraka_{\beta-\alpha}$.  It follows that the
ideal $\fraka_{\beta-\alpha}$ is prime.
\end{proof}

\begin{rem}
We have even proved that
$\fraka_{\beta-\alpha}=\fraka_\beta:\ideal{\del^\boldu}$ for a suitable
$\boldu$ is an associated prime of $\fraka_{\beta}$.
\end{rem}

\begin{ex}
Let $A=(1)$ and $\beta\in\CC$. Considering $H_A(\beta)$ for
$\beta\in\NN$ one easily
checks that $\fraka_\beta=\ideal{\del_1^{\beta+1}}$. While $\del_1:H_A(0)\to
H_A(1)\to H_A(2)\to\cdots$ are all isomorphisms, $\del_1:H_A(-1)\to
H_A(0)$ is not since $\fraka_0=\ideal{\del_1}$.
\end{ex}

We now introduce notation for each subset $F$ of $A$.  
\begin{df} Let
$F$ be a subset of $A$.  We write $\del_F=\{\del_k\st \bolda_k\in F\}$, $R_F$
for $\CC[\del_F]$, put $x_F=\{x_k\st \bolda_k\in F\}$, $\calO_F=\CC[x_F]$ and
denote by $D_F$ the ring of $\CC$-linear differential operators on
$\calO_F$. The {\em lattice of $F$} is the group $\ZZ F$, and we let
$E_F$ stand for the truncated Euler operators $\{\sum_{\bolda_k\in F}
a_{i,k} x_k\del_k\}_{i=1}^d$.

If $F$ is the set of all columns $\bolda\in A$ such that $L(\bolda)=0$
for some fixed functional $L:\ZZ^d\to \ZZ$ then we say that $F$ is a
{\em prime face} and define the ideals $I_F=\ker(R_F\to \CC[\NN
F])\subseteq R_F$, and $I^F_A=\ideal{\del_k\st \bolda_k\not\in
F}+R_A\cdot I_F\subseteq R_A$. Note the natural isomorphism
$R_A/I^F_A\cong R_F/I_F$. The collection of all ideals $I^F_A$ (where
$F$ runs through all prime faces of $A$) is identical with the
collection of all $\ZZ A$-graded prime ideals of $R_A$ that contain
$I_A$, which is in turn identified with the faces of the cone $\RR_+
A$.
\end{df}

\begin{df}
We say that $A$ is a {\em pyramid} over the prime face $F$ if
 $|A\setminus F|=1$.  In such a case, $\rk(F)=\rk(A)-1$.  We call $A$
 a {\em $t$-fold iterated pyramid} if there is a prime face $F$ of $A$
 such that $|A\setminus F|=k$ and $\rk(A)=\rk(F)+k$. This is
 equivalent to the existence of a sequence of prime faces
 $A=F_t\supseteq \ldots\supseteq F_0=F$ such that $F_{i+1}$ is a
 pyramid over $F_i$ for all $i$.
\end{df}

We now return to our hypothetical irreducible rank-jumping module
$H_A(\beta)(x)$ of minimal column number. By Proposition
\ref{prop-prime-mu}, $H_A(\beta)(x)=\calH_0(\beta;S_A)(x)\cong
\calH_0(\beta-\alpha;S_A)(x)$ where $\fraka_{\beta-\alpha}$ is prime and the
isomorphism is given by right multiplication of
$\calH_0(\beta-\alpha;S_A)$ by $\del^\boldu$ with $A\cdot \boldu =\alpha$. 

Since $\fraka_{\beta-\alpha}$ is a prime $\ZZ A$-graded ideal in $R_A$
then $\fraka_{\beta-\alpha}=I^F_A$ for some face $F$ obtained from $A$
by erasing the columns whose variables are linear generators of
$\fraka_{\beta-\alpha}$.  In consequence,
\begin{eqnarray*}
H_A(\beta-\alpha)(x)&=& D_A(x)/\ideal{I_A,E-\beta+\alpha}\\
&=&D_A(x)/\ideal{I^F_A,E-\beta+\alpha}\qquad 
\text{ since $I^F_A=\fraka_{\beta-\alpha}$}\\
&\cong&D_F(x_F)/\ideal{I_F,E_F-\beta+\alpha}\otimes_\CC D_{A\setminus
F}(x_{A\setminus F})/\ideal{\{\del_k\st \bolda_k\not\in F\}}\\ 
&=&H_F(\beta-\alpha)(x_F)\otimes_\CC \CC(x_{A\setminus F}). 
\end{eqnarray*}
We note that the second factor is irreducible of rank one and consider
now the hypergeometric system associated to $F$ and $\beta-\alpha$.
\begin{lem}
\mylabel{lem-vol-of-subset} Let $F$ be a subset of $A$. The simplicial
volume of $F$ cannot exceed the simplicial volume of $A$ (in their
respective lattices). Moreover, suppose that $F$ is a prime face of
$A$ such that the simplicial volume of $A$ and the simplicial volume
of $F$ (in their respective lattices) are the same. Then $A$ is an
iterated pyramid over $F$.
\end{lem}
\begin{proof}
Consider the following procedure.  Let $F_0=F$, let $L_0$ be the
lattice spanned by $F$ and let $V_0$ be the simplicial volume
of $F_0$.  Let $\bolda_{j_1},\ldots,\bolda_{j_t}$ be the columns of
$A\setminus F$.

For each triple $(F_i,L_i,V_i)$ already constructed, put
$F_{i+1}=F_i\cup\{\bolda_{j_i}\}$, and let $L_{i+1}$ and $V_{i+1}$ be
lattice and simplicial volume of $F_{i+1}$.

If $\dim(F_i)<\dim(F_{i+1})$ then $L_{i+1}\cong \ZZ\cdot
\bolda_{j_{i+1}}\oplus L_i$. Hence $V_i=V_{i+1}$. On the other hand, if
$\dim(F_i)=\dim(F_{i+1})$ then either $L_i=L_{i+1}$, or $L_{i+1}$
contains $L_i$ as a finite index subgroup. In the former case the
simplicial volumes $V_i$ and $V_{i+1}$ are computed with reference to
the same lattice, and since $F_{i+1}$ contains $F_i$ we see that
$V_i\le V_{i+1}$. In the second case even $[L_{i+1}:L_i]\cdot V_i\le
V_{i+1}$ follows.  Hence $V_0\le V_1\le \cdots \le V_t$ and the first part
of the lemma follows.

If now $F=F_t$ is a prime face 
with $\vol(F)=\vol(A)$ then $V_n=V_t$ and the equality
must hold throughout.  We claim that this proves the lemma. For,
assume that there is $i\geq t$ such that $F_{i+1}$ is not a pyramid
over $F_i$, and pick the smallest value $i\geq t$ of that
nature. Since the volumes are equal, we must have $\dim(F_i\cup
\{0\})=\dim(F_{i+1}\cup\{0\})$ as otherwise $F_{i+1}$ would be a
pyramid over $F_i$. 
Moreover,  we must have
$\bolda_{j_{i+1}}\in L_i\cap\conv(F_i\cup\{0\})$.  However, since
$F_i$ is an iterated pyramid over $F$ then the only points in the
intersection of $L_i$ with the convex hull of $F_{i}$ and the origin
are contained in $L_t\cup\{\bolda_{j_1},\ldots,\bolda_{j_i}\}$. Since
$F$ is a prime face, $\bolda_{j_{i+1}}$ cannot lie in $L_t$, and the
lemma follows.
\end{proof}

Thus, returning to $H_F(\beta-\alpha)$ in the theorem, 
\begin{eqnarray*}
\rk(H_{F}(\beta-\alpha))&=&\dim_{\CC(x_{F})}(H_{F}(\beta-\alpha)(x_{F}))\\
&=&\dim_{\CC(x_{A
    })}(\CC(x_A)
\otimes_{\CC[x_{F}]} H_{F}(\beta-\alpha))\\
&=&\dim_{\CC(x_{A
    })}(H_{A}(\beta)(x_A))\\
&=&\rk(H_A(\beta))>\vol(A)\geq \vol(F).
\end{eqnarray*}
In order to use the minimality condition of our supposed
counterexample we need to modify $(F,\beta-\alpha)$ to a case of full
rank.  Let $(A',\beta')$ be constructed as follows. First apply a
generic element $g$ of $GL(d,\ZZ)$ to the matrix
$(A,\beta-\alpha)$. The effect is that now the $\rk(F)$ top rows $A'$
of $g(F)$ are linearly independent, and the other $\rk(A)-\rk(F)$ rows
are $\ZZ$-linear combinations of the rows of $A'$. The latter implies
that the simplicial volumes of $A'$ and of $F$ are identical
(although, of course, their convex hulls live in spaces of different
dimension). Put $\beta'$ to be the top $\rk(F)$ rows of
$g(\beta-\alpha)$. Then the rows of $(F,\beta-\alpha)$ are in the
rowspan of $(A',\beta')$ and so $H_{A'}(\beta')=H_{F}(\beta-\alpha)$.

Above we have shown that
$\rk(H_{A'}(\beta'))=\rk(H_F(\beta-\alpha))>\vol(F)=\vol(A')$.  Hence
$H_{A'}(\beta')$ is rank-jumping with fewer columns than $A$. It
follows that $(A',\beta')$ is not a counterexample to the theorem and
hence $H_{A'}(\beta')(x_{F})$ has a nontrivial submodule. Any
associated nontrivial sequence
\[
0\to M\to H_{A'}(\beta')(x_{F})\to N\to 0,
\]
when tensored over $\CC$ with $\CC(x_{A\setminus F})$, gives a
nontrivial submodule for $H_A(\beta-\alpha)(x)\cong H_A(\beta)(x)$,
contradicting our hypothesis that $H_A(\beta)(x)$ be irreducible.
\end{proof}

\subsection{Reducibility for varying $\beta$}
In this subsection we consider to what extent the property of
$H_A(\beta)$ being irreducible is preserved under variation of the
parameter. 
\begin{cor}[to the proof of Theorem \ref{thm-jump-red}]
\mylabel{cor-irred-pyramid} 
Assume that $H_A(\beta)$ is irreducible.
Then there is a face $F$ of $A$ (this may be the trivial face $F=A$)
and $\boldu\in\NN^n$ with $A\cdot \boldu=\alpha\in\NN A$ such that
\begin{enumerate}
\item $A$ is an iterated pyramid over $F$,
\item $\del^\boldu: H_A(\beta-\alpha)(x)\to H_A(\beta)(x)$ is an
      isomorphism,
\item $\fraka_{\beta-\alpha}=I^F_A$ and $\{P\in R_F\st  P=0 \text{ in }
      H_A(\beta-\alpha)(x)\}=I_F$.
\end{enumerate}
\end{cor}
\begin{proof}
Suppose $H_A(\beta)$ is irreducible and $\fraka_{\beta}\not =0$. In
Proposition \ref{prop-prime-mu} we have proved that there is $\alpha\in
\NN A$ such that $H_A(\beta)\cong H_A(\beta-\alpha)$, the isomorphism is
given by right multiplication on $H_A(\beta-\alpha)$ by $\del^\boldu$ with
$A\cdot \boldu =\alpha$, and $\fraka_{\beta-\alpha}=I^F_A$ is the prime
ideal attached to the face $F$.
  
Let $\bar F=A\setminus F$. Our hypotheses imply that 
\begin{eqnarray*}
H_A(\beta-\alpha)(x)&=&D_A/D_A(I_A,E-\beta+\alpha)(x)\\
&=&D_A/D_A(I_A^F,E-\beta+\alpha)(x)\\
&=&D_F/D_F(I_F,E_F-\beta+\alpha)(x_F)\otimes_\CC \CC(x_{\bar F}).
\end{eqnarray*}
As at the end of the proof of Theorem \ref{thm-jump-red}, left multiply
$(A,\beta-\alpha)$ by a generic element of $GL(d,\ZZ)$. Then the top
$\dim(F)$ rows $(A',\beta')$ of $(gF,g\beta-g\alpha)$ give a
GKZ-system $H_{A'}(\beta')$ that is isomorphic to $H_F(\beta-\alpha)$
and where $A'$ has full rank.  In particular,
$H_{A'}(\beta')(x_F)\otimes_\CC \CC(x_{\bar F})\cong H_A(\beta)$ and
$\rk(H_{A'}(\beta'))=\rk(H_A(\beta))$. Now $H_A(\beta)$, and hence
$H_{A'}(\beta')(x_F)$, is irreducible and therefore $(A,\beta)$ and
$(A',\beta')$ are not rank jumps by Theorem \ref{thm-jump-red}. So
$\vol(F)=\vol(A')=\rk(H_{A'}(\beta'))=\rk(H_A(\beta))=\vol(A)$.

By Lemma \ref{lem-vol-of-subset}, $A$ is thus an iterated pyramid over
$F$. Finally, by construction $\fraka_{\beta-\alpha}=I^F_A$ and
$I^F_A\cap R_F=I_F$.
\end{proof}

We now state the main theorem in this subsection. 

\begin{thm}\label{thm-reducible-lattice}
If $H_A(\beta)$ is irreducible then $H_A(\gamma)$ is irreducible for
all $\gamma\in\beta+\ZZ^d$.
\end{thm}
\begin{proof}
Let $F$ be the smallest prime face of $A$ such that $A$ is an iterated
pyramid over $F$. Let $t=\dim(A)-\dim(F)$. Suppose that the columns of
$F$ are the final columns of $A$. By means of an element $g\in
GL(\ZZ^d)$ we can transform $(A,\beta)$ into a matrix with block
decomposition
$\left(\begin{array}{ccc}\id_t&F''&\beta''\\0&F'&\beta'\end{array}\right)$
where $\id_t$ is the $t\times t$ identity matrix.  Since $A$ is an
iterated pyramid over $F$, $gA$ is an iterated pyramid over $gF$. In
particular, the rank of $F'$ and the rank of
$\left(\begin{array}{cc}F''\\F'\end{array}\right)$ are equal. This
means that $F''$ is in the rowspan of $F'$ and so $g$ can be chosen in
such a way that $F''=0$. With this particular $g$ we consider
$H_{gA}(g\beta)$. Since $H_A(\beta)\cong H_{gA}(g\beta)$ it is
sufficient to prove the theorem under the assumption that $g=1$ and
$A$ is already in the block form
$\left(\begin{array}{cc}\id_t&0\\0&F\end{array}\right)$.

There is an isomorphism
$H_A(\beta)=D_t/D_t\cdot(\{x_i\del_i-\beta_i\}_{i=1}^t)\otimes_\CC
H_F(\beta_F)$ where $D_t$ is the Weyl algebra in $x_1,\ldots,x_t$, and
where $\beta_F$ are the last $\dim(F)$ rows of $\beta$. The module
$D_t/D_t\cdot(\{x_i\del_i-\beta_i\}_{i=1}^t)(x_F)$ is always irreducible,
so all reducibility issues of $H_A(\beta)$ are determined by those of
$H_F(\beta_F)$. According to our assumption, $H_F(\beta_F)$ is
irreducible, and we wish to show that $H_F(\gamma_F)$ is irreducible
for all $\gamma\in\beta+\ZZ F$. This reduces the proof of the theorem
to the case that $A$ is not a pyramid at all.

\bigskip

We now assume that $A$ is not a pyramid and that $H_A(\beta)$ is
irreducible. Let $\gamma\in\beta+\ZZ^d$; we shall show that
$H_A(\gamma)$ is irreducible as well.  By Corollary
\ref{cor-irred-pyramid}, $A$ is a pyramid unless
$\fraka_{\beta}=I_A$. We conclude that $\fraka_\beta=I_A$ and so all
maps $H_A(\beta-\alpha)\to H_A(\beta)$ given by right multiplication
by $\del^\boldu$ with $\alpha=A\cdot\boldu$ are isomorphisms.  Thus
$H_A(\beta)$ is isomorphic to $H_A(\gamma)$ if $\gamma\in \beta-\NN
A$.

It is hence sufficient to show that $H_A(\gamma')$ is irreducible
for some $\gamma'\in \gamma+\NN A$. Consider the intersection of
$\beta+\NN A$ and $\gamma+\NN A$. Since $\NN A$ generates $\ZZ^d$ as
a group, this intersection will contain a subset of the form
$\gamma''+\NN A$. Since the set of exceptional parameters is a Zariski
closed subset of positive codimension in $\CC^d$, there is $\gamma'\in
(\gamma+\NN A)\cap (\beta+\NN A)$ such that $\gamma'$ is not a rank
jump. Let $\boldu$ be an element of $\NN^n$ such that
$\beta+A\cdot\boldu=\gamma'$ and put $A\cdot \boldu=\alpha$.

Consider the short exact sequence  
\[
0\to S_A(\alpha)\to S_A\to S_A/\ideal{\del^\boldu}\to 0
\]
and the following part of the long exact Euler--Koszul sequence with
parameter $\gamma'$:
\begin{eqnarray}
\label{eqn-4-term}~\\\nonumber
\calH_1(\gamma';S_A/\ideal{\del^\boldu})\to
\calH_0(\gamma'-\alpha;S_A)\stackrel{\phi}{\longrightarrow} 
\calH_0(\gamma';S_A)\to \calH_0(\gamma';S_A/\ideal{\del^\boldu})\to 0.
\end{eqnarray}
Since $H_A(\gamma'-\alpha)(x)=H_A(\beta)(x)$ is irreducible, the map
$\phi(x)$ is either zero or injective.

We now prove that $\phi(x)$ cannot be zero. For, assume
otherwise. Then $H_A(\beta)(x)$ is a quotient of
$\calH_1(\gamma';S_A/\ideal{\del^\boldu})(x)$. Consider the
Euler--Koszul complex on the module $S_A/\ideal{\del^\boldu}$ and let
$e_1,\ldots,e_d$ be the canonical generators for the module
$\calK_1(\gamma';S_A/\ideal{\del^\boldu})$. As $\del^\boldu e_k=0$,
for all $k$ and for all $P=x^\boldv\del^\boldw \in D_A$ one has
$\del^\boldu Pe_k=[\del^\boldu,P]e_k$. Of course, the $x$-degree of
$[\del^\boldu,P]$ is less than that of $x^\boldv$. By iteration one
sees that the coset of $Pe_k$ is $\del^\boldu$-torsion. It follows
that $\calH_1(\gamma';S_A/\ideal{\del^\boldu})(x)$, and hence its
image in $H_A(\beta)(x)$, is $\del^\boldu$-torsion.  If $\phi(x)$ is
zero, the coset of $1$ in $H_A(\beta)$ is annihilated by the product
$f\del^{k\boldu}$ of a nonzero polynomial $f=f(x_1,\ldots,x_n)$ and a
high power of $\del^\boldu$. Since $H_A(\beta)(x)$ is irreducible and
$A$ is not a pyramid, the coset of $\del^{k\boldu}$ generates
$H_A(\beta)(x)$ for all $k\in\NN$. Since $f\del^{k\boldu}$ is zero in
$H_A(\beta)$, the support of $H_A(\beta)$ is then contained in the
variety of $f$. In that case, its rank must be zero since rank is
measured in a generic point, hence not in a point of $\var(f)$. But
$\rk(H_A(\beta))\geq\vol(A\cup\{0\})\geq \vol(\{0\})=1$. This
contradiction shows that $\phi(x)$ is nonzero and hence injective.

Since $\phi(x)$ is injective, the rank of $H_A(\beta)$ is a lower
bound for the rank of $H_A(\gamma')$. However, $H_A(\gamma')$ is not
rank-jumping and so $\rk(H_A(\beta))=\rk(H_A(\gamma'))$. By the
sequence (\ref{eqn-4-term}) the rank of
$\calH_0(\gamma';S_A/\ideal{\del^\boldu})$ is zero which by
Proposition 5.3 of \cite{MMW} means that
$\calH_i(\gamma';S_A/\ideal{\del^\boldu})=0$ for all $i$. The long
exact sequence (\ref{eqn-4-term}) gives now that $H_A(\beta)(x)\cong
H_A(\gamma')(x)$ is irreducible.  This concludes the proof.
\end{proof}
\begin{rem}
It would be interesting to know how irreducibility is encoded in
M.~Saito's invariants $E_\tau(\beta)$ discussed in
\cite{Saito-Compositio}. 
\end{rem}
\begin{rem}
Mellin systems are special systems of differential equations that are
satisfied by the roots of a polynomial $f$ with generic coefficients. 
In \cite{Dickenstein-Sadykov}, A.~Dickenstein and T.~Sadykov prove by
ingeneous methods that every Mellin system has reducible monodromy,
i.e., is not the minimal differential system for the roots. 

As it turns out, every Mellin system surjects onto an
$A$-hypergeometric system $H_A(\beta)$ where the parameter $\beta$ has
only {\em integral} coordinates, determined by the support of $f$. By
Theorem \ref{thm-reducible-lattice}, the reducibility of any such
$H_A(\beta)$, and hence the reducibility of every Mellin system can be
reduced to the case $\beta=0$.
\begin{prop}
Let $A$ be such that it is not a pyramid over a point. (In other
words, assume that $\vol(A)>1$.)
The system $H_A(\beta)$ is reducible provided that $\beta=0$.
\end{prop}
\begin{proof}
If $\beta=0$ then $I_A$ and $E-\beta$ are both contained in the ideal
$\del_A=\ideal{\del_1,\ldots,\del_n}$. Hence there is a surjection 
\[
0\to K\to H_A(\beta)\onto D_A/\del_A\to 0
\]
with kernel $K$.
Now rank is additive in short exact sequences. Since
$\rk(H_A(\beta))=\vol(A)>1=\rk(D_A/\del_A)$, the rank of $K$ must be
positive; hence $H_A(\beta)$ has reducible monodromy.
\end{proof}
We note that Mellin systems do not lead to $\vol(A)=1$ apart from the
completely degenerate case of a monomial $f$. 
\end{rem}

\section{Duality for hypergeometric systems}

In this section we investigate holonomic duality of GKZ-systems. 
We show that in the Gorenstein case holonomic duality is a functor that
preserves the set of GKZ-systems. 
More generally, holonomic duals of
GKZ-systems are GKZ-systems for all $A$ and generic $\beta$,
but for special $A,\beta$ this is false. 

We begin with the simplest possible case, that of a Gorenstein toric
variety.  Recall that $\eps_A=\sum_{j=1}^n \bolda_j$ is the canonical
degree of $R_A$. 
\begin{prop}
Suppose that $A$ is Gorenstein. Then $\DD(H_A(\beta))\cong
H_A(-\beta-\eps_A-c_A)$ is a proper GKZ-system for all $\beta$, where
$c_A$ is the canonical degree of $S_A$.
\end{prop}
\begin{proof}
We consider the duality spectral sequence (\ref{eqn-ss}).  Since $A$
is Cohen--Macaulay, we obtain duality
\[
\DD(H_A(\beta))^\vee=\calH_0(-\beta-\eps_A,\ext^{n-d}_{R_A}(S_A,R_A)).
\] 
Since $S_A$ is Gorenstein, $\ext^{n-d}_{R_A}(S_A,R_A)$ is isomorphic to
$S_A(c_A)$ where $c_A$ is the canonical degree of $S_A$. Then
$\DD(H_A(\beta))^\vee \cong H_A(-\beta-c_A-\eps_A)$.
\end{proof}
It is clear that any time the duality spectral sequence is used, we
obtain an identification of the dual of a GKZ-system with another
GKZ-system only up to a coordinate transformation. In the projective
case, this transformation can be omitted.
\begin{ex}
Consider the curve defined by $A=(a_1,\ldots,a_n)$ with $0\le a_i\le
a_{i+1}$ for all $i$. If $S_A$ is Gorenstein and if $c$ is the largest
integer in $\NN\setminus \NN A$ then $\ext^{n-1}_{R_A}(S_A,R_A)$ is
isomorphic to $S_A(c)$. It follows that $\DD(H_A(\beta))^\vee\cong
\calH_0(-\beta-\eps_A;S_A(c))= H_A(-\beta-\eps_A-c)$.

The shift that occurs (besides negation of $\beta$) depends not only
on the obvious parameters $n$, $d$ and $\eps_A$, but also on the
arithmetic of the semigroup $\NN A$. For example, $A=(2,5)$ and
$A=(3,4)$ have shifts by $-7-3$ and $-7-5$ respectively. 
Both are complete intersections and the shift is precisely
minus the (homogenized) degree of the curve described ($2\times 5$ and
$3\times 4$ respectively).  If $A$ is a complete intersection on more
than two powers of $x$ then the shift is the least common multiple of
the (homogenized) degrees of the defining equations.
\end{ex} 

\medskip

To illustrate the difficulties that one struggles with when testing
whether a given module is isomorphic to a GKZ-system, let us consider
the following

\begin{ex}
Let $A=\left(\begin{array}{cccc}1&1&1&1\\0&1&2&3\end{array}\right)$,
and $\beta\in\CC^2$; then
\[
I_A=\ideal{\del_1\del_3-\del_2^2,\del_2\del_4-\del_3^2,
\del_1\del_4-\del_2\del_3}
\]
is the ideal of the projective rational cubic, which is arithmetically
Cohen--Macaulay. By the spectral sequence (\ref{eqn-ss}) above,
$\DD(H_A(\beta))$ is $D_A\otimes_{R_A}\ext^2_{R_A}(S_A,R_A)$ modulo
the images of the endomorphisms
\[
E_1-\beta_1:m\mapsto
-(4+x_1\del_1+x_2\del_2+x_3\del_3+x_4\del_4+\beta_1+\deg_1(m))m
\]
and
\[
E_2-\beta_2:m\mapsto
-(6+x_2\del_2+2x_3\del_3+3x_4\del_4+\beta_2+\deg_2(m))m.
\]  
In particular, $\DD(H_A(\beta))$ is a generalized $A$-hypergeometric
system.

By explicit computation one sees with $M=\ext^3_{R_A}(S_A,R_A)$ that
\[
M\cong
D_A^2/\ideal{(\del_4,-\del_3),(\del_3,-\del_2),(\del_2,-\del_1)}
\]
where $m_1,m_2$ are the cosets of  $(1,0)$ and
$(0,1)$ with degrees $(-3,-5)$ and $(-3,-4)$
respectively.

Since $\ann_{D_A}(m_1)=D_A\cdot I_A$, it follows
that $M':=M/D_A\cdot m_1=D_A/\ideal{\del_3,\del_2,\del_1}$. Hence with
$\bar m_2=m_2
\mod D_A\cdot m_1$ we see that
\begin{eqnarray*}
(E_2-3E_1)\circ \bar m_2&=&
(3x_1\del_1+2x_2\del_2+x_3\del_3+3\deg_1(m_2)+3\beta_1-\deg_2(m_2)-\beta_2)
\bar m_2\\
&=&(3\beta_1-\beta_2+3\deg_1(m_2)-\deg_2(m_2))\bar m_2.
\end{eqnarray*}
Therefore if $\beta_2-4\not =3\beta_1-3\cdot 3$ then
$3E_1-E_2$ is a unit on $M'$ and so $M/\im(E-\beta)=D_A\cdot
m_1/(D_A(E_1-\beta_1)m_1+D_A(E_2-\beta_2)m_1)\cong
H_A((4-3,6-5))$ is a proper $A$-hypergeometric system.

On the other hand, put $M'':=M/D_A\cdot
m_2=D_A/\ideal{\del_4,\del_3,\del_2}$. Hence the $M$-endomorphism
$E_2-\beta_1$ is multiplication by $-\beta_2-\deg_2(m_1)$ on $M''$ and
as long as $-\beta_2+5\not =0$ we conclude that
$M/\im(E_1-\beta_1,E_2-\beta_2)=D_A\cdot m_2/((E_1-\beta_1)\circ D_A
m_2 +(E_2-\beta_2)\circ D_A m_2)\cong H_A((4-3,6-5))$ is a
proper $A$-hypergeometric system as well.

Let finally $\beta_2-5=0$, $\beta_2-4 =3\beta_1-3\cdot 3$, so
$\beta=(10/3,5)$. We view $N=\ext^2_{R_A}(S_A,R_A)$ as the ideal in
$S_A$ generated by $\del_2$ and $\del_3$, shifted by $(-2,-3)$ so that
$\deg(\del_2)=(-3,-4)$ and $\deg(\del_3)=(-3,-5)$. Then
\[
S_A(2,3)/N=R_A/\ideal{\del_1\del_4,\del_2,\del_3}(2,3).
\]
Since no element of $R_A/\ideal{\del_2,\del_3,\del_4}$ has fractional
degree, the endomorphism $E_2$ is a unit on
$D_A/\ideal{\del_2,\del_3,\del_4}(2,3)$.  At the same time, $E_2-3E_1$
is a unit on $D_A/\ideal{\del_1,\del_2,\del_3}(2,3)$ since the degrees
in $R_A/\ideal{\del_1,\del_2,\del_3}(2,3)$ are of the form
$(a-2,3a-3)$ while $\beta=(10/3,5)$ is not. So
$M/\im(E_1-\beta_1,E_2-\beta_2)\cong
S_A(2,3)/\im(E_1-\beta_1,E_2-\beta_2)=H_A(-\beta_1,-\beta_2)=H_A(-16/3,-8)$.
\end{ex}

The example shows that even if $\ext^d_{R_A}(S_A,R_A)$ is not cyclic
(as it is in the Gorenstein case), the dual of $H_A(\beta)$ may often
be identified with a hypergeometric system. But there is an obvious
difficulty in picking the right isomorphism as there does not seem to
be a $\beta$-global one for non-Gorenstein $A$.

We consider next the case where $S_A$ is normal. Let $\{H_1\geq
0,\ldots,H_e\geq 0\}$ be the defining hyperplanes of $\RR_{\geq
0}A\subseteq \RR^d$.
\begin{prop}
Suppose $\NN A$ is normal. Then the set of all $\gamma$ for which
$\DD(H_A(\gamma))$ is not a proper GKZ-system is contained in a finite
subspace arrangement of codimension at least three. All participating
subspaces are shifts of faces of $\RR_{\geq 0}A$ by lattice elements.
\end{prop}
\begin{proof}
If $A$ is normal then it is Cohen--Macaulay and hence the duality
spectral sequence collapses: $\DD(H_A(\gamma))\cong
\calH_0(-\gamma-\eps_A;\ext^{n-d}_{R_A}(S_A,R_A))$ up to a coordinate
change.  Moreover, the canonical module $\ext^{n-d}_{R_A}(S_A,R_A)$ is
the interior ideal $N_A$ of $S_A$, with an appropriate shift. We
therefore want to show that for all choices $\beta$ away from a
codimension three arrangement the module $\calH_0(\beta;N_A)$ is a
proper GKZ-system.

Consider the natural inclusion $N_A\into S_A$ and the quotient
$Q$. This is a toric sequence and hence there is a long exact
sequence
\[
\cdots\to \calH_{1}(\beta;Q)\to \calH_0(\beta;N_A)\to \calH_0(\beta;S_A)\to
\calH_0(\beta;Q)\to 0.
\]
The quasi-degrees of $Q$ are those elements of $\CC^d$ that sit on a
complexified defining hyperplane of $\RR_{\geq 0} A$.  Hence if
$-\beta$ is not on any complexified defining hyperplane of $\RR_{\geq
0}A$, then $\calH_0(\beta;N_A)$ and $H_A(\beta)=\calH_0(\beta;S_A)$
are isomorphic.

Let now $-\beta\not =0$ sit on one such hyperplane. Then there is a
smallest infinite face $F$ such that $-\beta$ is a quasi-degree of $F$
but not of any smaller face. Note that 
there is an infinite face $N_F$ of $\conv(\deg(N_A))$  such that some shift 
 $F+z_F$ of $F$, $z\in\ZZ^d$,
is inside $N_F$ while $\dim(N_F/F+z_F)<\dim(N_F)=\dim(F)$.
Let $f_F$ be a point in $\tdeg(N_A)$ sitting on $F+z_F$. Let
$\delta_F$ be an element of $\NN A$ in the relative interior of $F$
itself. Then $f_F-l_F\delta_F+\NN A$ will contain $\tdeg(N_A)$ for
some large $l_F\in\NN$. We hence have an injection $N_A\into
S_A(l_F\delta_F-f_F)$ and we let $Q_F$ be the cokernel. Now note that
since $t^{l_F\delta_F}\cdot S_A(l_F\delta_F-f_F)\subseteq S_A$, all
components of $\qdeg (Q_F)$ are lattice shifts of $\CC F'$ where $F'$
does not contain $F$.

Therefore, the quasi-degrees of $Q_F$ are
contained in a finite union of hyperplanes of the form $H_j=r_{j,k,F}$
where $H_j$ is a defining hyperplane of $\NN A$ transversal to $F$ and
$H_j(f_F-l_F\delta_F)\le r_{j,k,F}<\min\{r\in H_j(N_A)\}$.
Note that this construction depends only on the minimal face $F$
containing $-\beta$, so there are finitely many $r_{j,k,F}$. 

Now consider the embedding $S_A(-f_F)\into N_A$. The true degrees, and
hence the quasi-degrees, of the cokernel are obviously contained in a
union of hyperplanes of the form $H_j=s_{j,k,F}$ where
$H_j(f_F)>s_{j,k,F}\geq\min\{s\in H_j(N_A)\}$ and $H_j$ is
perpendicular to $F$.

We conclude that for $\beta\not =0$ the module $H_0(\beta;N_A)$ is
isomorphic to a proper GKZ-system provided that
\begin{itemize}
\item $\beta$ is not on $H_j=0$ for any defining hypersurface of $\NN
  A$, or
\item for the minimal face $F$ with $-\beta\in\qdeg(F)$ and all $H_j$
  transversal to $F$ with $H_j(\beta)\in H_j(\ZZ A)$ we have
  $H_j(f_F-l_F\delta_F)\le H_j(\beta)<\min\{r\in H_j(N_A)\}$, or
\item for the minimal face $F$ with $-\beta\in\qdeg(F)$ and all $H_j$
  transversal to $F$ with $H_j(\beta)\in H_j(\ZZ A)$ we have
  $H_j(f_F)> H_j(\beta)\geq \min\{s\in H_j(N_A)\}$
\end{itemize}
because in all those cases $H_0(\beta;Q)=0$.
Hence, if $\beta\not =0$ then $\calH_0(\beta;N_A)$ is a proper
GKZ-system unless $\beta$ violates all three (disjoint) conditions
listed above.  It remains to show that $\calH_0(0;N_A)$ is a
GKZ-system. Pick an element $\boldu\in N_A$ and consider $0\to
S_A(-\deg(\del^\boldu))\into N_A\onto Q\to 0$ for a suitable toric
module $Q$. Suppose $0$ is a quasi-degree of $Q$. As $0\not\in N_A$,
there is a face $F$ of $\NN A$ such that $\RR F$ meets
$\deg(N_A)$. Obviously the only face that is capable of doing this is
the whole semigroup $\NN A$. On the other hand, $Q$ is not
$d$-dimensional since $N_A/S_A\cdot\del^\boldu\subseteq S_A/S_A\cdot
\del^\boldu$ is $(d-1)$-dimensional. Hence $\beta=0$ is not a
quasi-degree of $Q$ and again the long exact sequence tells us that
$\calH_0(0;N_A)\cong H_A(\deg(\del^\boldu))$.

It follows that the set of parameters for which $\DD(H_A(\beta))$ is
not a proper GKZ-system is contained in an arrangement of codimension
three, and the origin is not part of this set.
\end{proof}

\begin{cor}
If $\dim(S_A)\le 2$ and $S_A$ is normal then the holonomic duals of all
$H_A(\beta)$ are proper GKZ-systems.  \qed
\end{cor}
Normal semigroups are Cohen--Macaulay by Hochster's theorem. 
We now show by example that absence of 
Cohen--Macaulayness may bring with it the failure of Takayama's
conjecture.

\begin{prop}
Let
$A=\left(\begin{array}{cccc}1&1&1&1\\0&1&3&4\end{array}\right)$ and
$M=\DD(H_A((1,2)))$. Then $M$ is not a GKZ-system
in the classical sense. That is, $M\not = H_B(\beta')$ for
any $B$ and $\beta'$.
\begin{proof}
The only candidate for $B$ is $A$ itself, or a matrix that is obtained
from $A$ by an action of $GL(\ZZ,2)$. We may hence assume that
$B=A$. Denote by $X_i$ the module $\ext^i_{R_A}(S_A,R_A)$.  So $X_2$
is the canonical module of the normalization of $S_A$.  It is
well-known that $A$ is not Cohen--Macaulay, and that $X_3\cong \CC$,
generated in degree $(5,10)$.

By \cite{Sturmfels-Takayama}, 
$H_A((1,2))$ is of rank five and $(1,2)$ is the only
exceptional parameter. By duality,
$M$ is of rank five as well, so only $H_A((1,2))$ is a
candidate for equaling $M$ among all $A$-hypergeometric systems.

Since $\eps_A=(4,8)$ is the canonical degree of $R_A$, and $(5,10)$ is
the unique degree of $X_3$, $E+(1,2)+\eps_A=E+(5,10)$ is the zero
morphism on $D_A\otimes_{R_A}X_3$.  Thus,
$\calH_{i}(-\beta-\eps_A;X_3)\cong \bigwedge^i(\calO_A\oplus\calO_A)$.
In the spectral sequence (\ref{eqn-ss}) there are hence four nonzero
terms in the $E_2^{p,q}$-picture: $\calH_0(-\beta-\eps_A;X_2)$, and
$\calH_i(-\beta-\eps_A;X_3)$ for $0\le i\le 2$.  The sequence
converges to $M=M^\vee =\DD(H_A((1,2)))^\vee$. Let
$\beta=(1,2)$. Since there are only two nonzero columns in
(\ref{eqn-ss}), there is a ``long'' exact sequence
\[
0\to \calH_{2}(-\beta-\eps_A;X_3)\to \calH_{0}(-\beta-\eps_A;X_2)\to
M\to \calH_{1}(-\beta-\eps_A;X_3)\to 0.
\]
The display shows that $M\to\calO_A\oplus\calO_A\to 0$ is exact
and hence $\hom_{D_A}(M,\calO_A)$ is at
least a 2-dimensional vector space. However, by Proposition 3.4.11 in
\cite{SST}, $\calO_A$ contains at most one polynomial solution to any
$H_A(\beta')$.  This contradiction shows that $M$ is not a GKZ-system
proper.
\end{proof}
\end{prop}

\begin{rem}
Let $A$ be a $2\times n$ matrix of rank two. If $\NN A$ is not
Cohen--Macaulay, then $N=\ext^{n-1}_{R_A}(S_A,R_A)$ is a
finite-dimensional nonzero vector space. The exceptional set of $A$ is
the set of degrees $\beta\in\tdeg(N)-\eps_A$. The module
$X_{n-2}=\ext^{n-2}_{R_A}(S_A,R_A)$ is Cohen--Macaulay, so has no
higher Euler--Koszul homology.

For each such exceptional $\beta$ the duality spectral sequence has
the form as in the proposition. Namely, there are two nonzero columns,
the $(n-2)$-nd and the $(n-1)$-st. The $(n-1)$-st contains three
modules, $\calO_A$, $\calO_A^2$, $\calO_A$, and the $(n-2)$-nd
contains one module. The spectral sequence shows an epimorphism from
$\DD(H_A(\beta))$ onto $\calO^2_A$. In particular, $A$-exceptional
degrees give duals that are not $A$-hypergeometric systems.
\end{rem}

In the light of the proposition it would be interesting to know what
can be rescued if $A$ is not normal but Cohen--Macaulay. However, no
matter how bad $A$ is, generic parameters satisfy Takayama's
conjecture: 

\begin{thm}
For all $A\in\ZZ^{d,n}$ there is a Zariski open set $U\subseteq \CC^d$
of parameters such that $\DD(H_A(\beta))$ is a proper GKZ-system for
all $\beta\in U$.
\end{thm}
\begin{proof}
We note that $H_A(\beta)=\calH_0(\beta;S_A)$ is for generic $\beta$
isomorphic to $\calH_0(\beta;\tilde S_A)$ since $\tilde S_A/S_A$ is an
$R_A$-module of dimension at most $d-1$. Since $\tilde S_A$ is a
Cohen--Macaulay $S_A$-module, the dual of $\calH_0(\beta;\tilde S_A)$
equals $\calH_0(-\beta-\eps_A,\ext^{n-d}_{R_A}(\tilde S_A,R_A))$ up to
a coordinate transformation.

Let $\tilde A\supseteq A$ be a matrix composed of columns that generate the
saturation of the semigroup $\NN A$: $\NN \tilde A=\RR_{\geq
0}A\cap\ZZ A$.  Let $R_{\tilde A}\supseteq R_A$ be the polynomial ring
$\CC[\del_{\tilde A}]$ surjecting onto $\tilde S_A=S_{\tilde A}$.
Then $\ext^{n-d}_{R_A}(\tilde S_A,R_A))$ and $\ext^{\tilde
n-d}_{R_{\tilde A}}(\tilde S_A,R_{\tilde A})$ are isomorphic as
$R_A$-modules.  To see this, note first that the restriction from
$R_{\tilde A}$-modules to $R_A$-modules preserves the degree.  Then
local duality over $R_A$ shows that $\ext^{n-d}_{R_A}(\tilde
S_A,R_A)^*\cong H_{\frakm_A}^d(\tilde S_A)$ over $R_A$. Here $(-)^*$
denotes the vector space dual on the graded pieces, followed by the
shift by $\eps_A$. Since the radical of $\frakm_A \tilde S_A$ is
$\frakm_{\tilde A}$, $H_{\frakm_A}^d(\tilde S_A)=H_{\frakm_{\tilde
A}}^d(\tilde S_A)$ as $R_A$-modules. By local duality over $R_{\tilde
A}$, $H^d_{\frakm_{\tilde A}}(\tilde S_A)$ is isomorphic to
$\left(\ext^{\tilde n-d}_{R_{\tilde A}}(\tilde S_A,R_{\tilde
A})\right)^*$ shifted by $\eps_A-\eps_{\tilde A}$. Thus,
$\ext^{n-d}_{R_A}(\tilde S_A,R_A)\cong \ext^{\tilde n-d}_{R_{\tilde
A}}(\tilde S_A,R_{\tilde A})(\eps_{\tilde A}-\eps_A)$. So
$\DD(H_A(\beta))\cong \calH_0(-\beta-\eps_{\tilde A},\ext^{\tilde
n-d}_{R_{\tilde A}}(\tilde S_A,R_{\tilde A}))$ for generic $\beta$ (as
$D_A$-modules).

Since $\tilde A$ is normal, we have $\ext^{\tilde n-d}_{R_{\tilde
A}}(\tilde S_A,R_{\tilde A})\cong N_{\tilde A}$, the interior ideal of
$\tilde S_A$. This module injects into $S_A$ with a cokernel of
dimension at most $d-1$. Hence for generic $\beta$,
$\DD(H_A(\beta))\cong \calH_0(-\beta-\eps_{\tilde A},\ext^{\tilde
n-d}_{R_{\tilde A}}(\tilde S_A,R_{\tilde
A}))=\calH_0(-\beta-\eps_{\tilde A},N_{\tilde A})\cong
H_A(-\beta-\eps_{\tilde A})$.
\end{proof}

Based on our results and computer experiments we close with
\begin{conj}
The set of parameters $\beta$ for which the dual of $H_A(\beta)$ is a
classical GKZ-system agrees with the points outside the
$A$-exceptional arrangement defined as the union of the quasi-degrees
of all local cohomology modules $H^i_\frakm(S_A)$ with $i<\dim(S_A)$.
\end{conj}

\bibliographystyle{abbrv}
\bibliography{../../bib.bib}

\end{document}